\documentclass[12pt]{article}
\usepackage{amsthm}
\usepackage{amsfonts}
\usepackage{amsmath}
\begin{document}

% Your \newcommands below (if there are any):
\theoremstyle{plain}
\newtheorem{theorem}{Theorem}
\newtheorem{definition}[theorem]{Definition}
\newtheorem{corollary}[theorem]{Corollary}
\newtheorem{lemma}[theorem]{Lemma}

\theoremstyle{definition}
\newtheorem{example}[theorem]{Example}
\newtheorem{remark}[theorem]{Remark}
\makeatletter
\def\rdots{\mathinner{\mkern1mu\raise\p@
\vbox{\kern7\p@\hbox{.}}\mkern2mu
\raise4\p@\hbox{.}\mkern2mu\raise7\p@\hbox{.}\mkern1mu}}
\makeatother
\renewcommand{\arraystretch}{1.8}

\oddsidemargin 16.5mm
\evensidemargin 16.5mm

\thispagestyle{plain}

\vspace{5cc}
\begin{center}

{\large\bf  COMPOSITA AND ITS PROPERTIES
\rule{0mm}{6mm}\renewcommand{\thefootnote}{}%Enter at least one, but not more than 3 MSCs.
% First entered MSC will be a primary one, others (at most 2) will be secondary.
\footnotetext{

\rule{2.4mm}{0mm}Keywords and Phrases. Generating functions, Composition, Composita, Functional equations.
}}

\vspace{1cc}
{\large\it Vladimir V. Kruchinin, Dmitry V. Kruchinin }

\vspace{1cc}
\parbox{24cc}{{\small

In this paper we study the coefficients of the powers of an ordinary generating function and their properties.
A new class of functions based on compositions of an integer $n$ is introduced and is termed composita.
We present theorems about compositae  and operations with compositae. We obtain the compositae of polynomials,  trigonometric and hyperbolic functions. Using the notion of the composita we get the solution of the functional equation $B(x)=H(xB(x)^m)$, where $H(x),B(x)$ are generating functions, and $m\in \mathbb{N}$.

}}
\end{center}

\vspace{1cc}

\vspace{1.5cc}
\begin{center}
{\bf 1. INTRODUCTION}
\end{center}

The computations based on combinatorial objects are an important direction of research in enumerative combinatorics and related fields of mathematics. For example, ordered partitions of a finite set is used to derive the formula for a composition of exponential generating functions \cite{Stanley_v2}. 
Computations that use compositions of an integer $n$ are found in various problems: derivation of a convolution of convolutions \cite{ConcreteMath}, composition of ordinary generating functions \cite{Egor}, calculation of the $n$-th order derivatives of a composite function \cite{WJohnson}, generation of ordered root trees \cite{Kru2010}, etc.
However, there is no unified approach to solving problems based on compositions. 

In this work, we consider a unified approach to the above problems, using a special function termed a \textit{composita}. 
The notion of the composita is close to that of a Riordan array \cite{Shapiro,Sprugnoli}, but the composita characterizes only one function, and  potential polynomials for exponential generating functions \cite{Comtet}.

Most of all papers and books related to combinatorial problems and generating functions use coefficients of the powers of an ordinary generating function \cite{Stanley_v2,Egor,Comtet,Flajolet}. However, as an independent object of study this has not considered. So investigation of the coefficients of the powers of an ordinary generating function is very important.

\vspace{1.5cc}
\begin{center}
{\bf 2. COMPOSITA}
\end{center}

Now we introduce the definition of composita.
\begin{definition} The composita of the generating function $F(x)=\sum_{n>0}f(n)x^n$ is the function of two variables
\begin{equation}
%\label{Fnk0}F_{n,k}^{\Delta}=\sum_{\pi_k \in C_n}{f(\lambda_1)f(\lambda_2)\ldots f(\lambda_k)}.
\label{Fnk0}F^{\Delta}(n,k)=\sum_{\pi_k \in C_n}{f(\lambda_1)f(\lambda_2)\ldots f(\lambda_k)},
\end{equation}
where $C_n$ is a set of all compositions of an integer $n$, $\pi_k$ is the composition $n$ into $k$ parts such that $\sum_{i=1}^k\lambda_i=n$.
\end{definition}

It follows from the definition of a composita that it is defined for a generating function $F(x)$ for which $f(0)=0$. Let us consider a generating function $F(x)=\frac{x}{1-x}=\sum_{n>0}x^n$. On the strength of formula (\ref{Fnk0}), the composita of this function is
$$
 F^{\Delta}(n,k)={n-1 \choose k-1}.
$$
For all $n>0$ we have $f(n)=1$; therefore, the formula (\ref{Fnk0}) counts the number of compositions of $n$ into $k$ parts. 

Next we obtain a recurrent formula for the composita of a generating function.

\begin{theorem} \label{TheoremFnk}
For the composita $F^{\Delta}(n,k)$ of the generating function $F(x)=\sum\limits_{n>0}f(n)x^n$ the following recurrent relation holds true
\begin{equation}
\label{CompositaRecur}
F^{\Delta}(n,k)=
\begin{cases}
f(n), & \text{if $k=1$};\\
\sum\limits_{i=1}^{n-k+1}f(i)F^{\Delta}(n-i,k-1), &
 \text{if $k\leq n$}.
\end{cases}
\end{equation}
\end{theorem}

\begin{proof}  
The composition $\pi_k$ for $k=1$ is unique and is equal to $n$; from whence it follows that $F^{\Delta}(n,1)=f(n)$. Now for $k>1$ we group in the formula (\ref{Fnk0}) all products $f(\lambda_1)f(\lambda_2)\ldots f(\lambda_k)$ of the composition $\pi_k$ with equal $\lambda_1$. Let us take $f(\lambda_1)$ out of the brackets; we see that the sum of the products in the brackets is equal to $F^{\Delta}(n-\lambda_1,k-1)$. Then for all values of $\lambda_1$ we obtain
\begin{multline*}
F^{\Delta}(n,k)=f(1)F^{\Delta}(n-1,k-1)+f(2)F^{\Delta}(n-2,k-1)+\cdots \\ \cdots +f(i)F^{\Delta}(n-i,k-1)+\ldots +f(n-(n-k+1))F^{\Delta}(k-1,k-1).
\end{multline*}

Thus, the theorem is proved.
\end{proof}

It is obviously that
$$
F^{\Delta}(n,n)=f(1)F^{\Delta}(n-1,n-1)=f(1)^n.
$$
The formula (\ref{CompositaRecur}) allows  the conclusion that the composita is a characteristic of the generating function $F(x)$.

In tabular form, the composita is presented as a triangle as follows
$$
\begin{array}{ccccccccccc}
&&&&& F_{1,1}^{\Delta}\\
&&&& F_{2,1}^{\Delta} && F_{2,2}^{\Delta}\\
&&& F_{3,1}^{\Delta} && F_{3,2}^{\Delta} && F_{3,3}^{\Delta}\\
&& F_{4,1}^{\Delta} && F_{4,2}^{\Delta} && F_{4,3}^{\Delta} && F_{4,4}^{\Delta}\\
& \rdots && \vdots && \vdots && \vdots && \ddots\\
F_{n,1}^{\Delta} && F_{n,2}^{\Delta} && \ldots && \ldots && F_{n,n-1}^{\Delta} && F_{n,n}^{\Delta}\\
\end{array}
$$
or, since $F_{1,n}^{\Delta}=f(n)$, $F_{n,n}^{\Delta}={[f(1)]}^n$, as
$$
\begin{array}{ccccccccccc}
&&&&& f(1)\\
&&&& f(2) && f^2(1)\\
&&& f(3) && F_{3,2}^{\Delta} && f^3(1)\\
&& f(4) && F_{4,2}^{\Delta} && F_{4,3}^{\Delta} && f^4(1)\\
& \rdots && \vdots && \vdots && \vdots && \ddots\\
f(n) && F_{n,2}^{\Delta} && \ldots && \ldots && F_{n,n-1}^{\Delta} && f^n(1)\\
\end{array}
$$ 

Presented below are the first terms of the composita of the generating function $F(x)=\frac{x}{1-x}$ (it is the Pascal triangle)
$$
\begin{array}{ccccccccccc}
&&&&&1\\
&&&& 1 && 1\\
&&& 1 && 2 && 1\\
&& 1 && 3 && 3 && 1\\
& 1 && 4 && 6 && 4 && 1\\
1 && 5 && 10 && 10 && 5 && 1
\end{array}
$$

For the given generating function $F(x)=\sum_{n>0}{f(n)x^n}$ the composita $F^{\Delta}(n,k)$ always exists and is unique.

Next we consider a generating function of the composita. 
The generating function of the composita of $F(x)$ is equal to
\begin{equation}
\label{GenComp}
[F(x)]^k=\sum_{n\geq k}F^{\Delta}(n,k)x^n.
\end{equation} 
It follows from 
$$
[F(x)]^k=\sum_{n\geq  k}\sum_{\pi_k \in C_n}{f(\lambda_1)f(\lambda_2)\ldots f(\lambda_k)}x^n=\sum_{n\geq  k}F^{\Delta}(n,k)x^n. 
$$

The composita is the coefficients of the powers of an ordinary generating function
$$F^{\Delta}(n,k) := [z^{n}] F(x)^k.$$

\vspace{1.5cc}
\begin{center}
{\bf 3. OPERATIONS WITH COMPOSITAE}
\end{center}

The above result allows us to use generating functions for computation of compositae. In this section we introduce several theorems for computation of compositae.

\begin{theorem}
\label{Theorem_shift}
 Suppose $F(x)=\sum_{n>0} f(n)x^n$ is a generating function, $F^{\Delta}(n,k)$ is the composita of this generating function. Then for the generating function $A(x)=xF(x)$  the composita is equal to 
\begin{equation}
A^{\Delta}(n,k)=F^{\Delta}(n-k,k).
\end{equation}
\end{theorem}

\begin{proof}
Using (\ref{GenComp}), we get
$$
[A(x)]^k=[xF(x)]^k=x^k[F(x)]^k=\sum_{m\geq  k} F^{\Delta}(n,k)x^{m+k}.
$$
Substituting $n$ for $m+k$, we get the following expression
$$[A(x)]^k=\sum_{n\geq  2k} F^{\Delta}(n-k,k)x^n.
$$
Therefore,
$$
A^{\Delta}(n,k)=F^{\Delta}(n-k,k).
$$
\end{proof}

\begin{corollary}
\label{compB(n,k)}
Suppose $B(x)=\sum_{n\geq 0} b(n)x^n$ is a generating function  such that $[B(x)]^k=\sum_{n\geq 0}B(n,k)x^n$. Then the composita of the generating function $A(x)=xB(x)$ is equal to
\begin{equation}
A^{\Delta}(n,k)=B(n-k,k).
\end{equation}
\end{corollary}

\begin{corollary}
\label{coeffB(n,k)}
Suppose $A(x)=\sum_{n>0} a(n)x^n$ is a generating function, $A^{\Delta}(n,k)$ is the composita of this generating function. Then for the generating function $[B(x)]^k=[\frac{F(x)}{x}]^k=\sum_{n\geq 0}B(n,k)x^n$ such that $B(x)=\sum_{n\geq 0} b(n)x^n$    the expression of coefficients is equal to 
\begin{equation}
B(n,k)=A^{\Delta}(n+k,k).
\end{equation}
\end{corollary}

\begin{theorem}
\label{comp-b}
Suppose $B(x)=\sum_{n\geq 0} b(n)x^n$ is a generating function  such that $[B(x)]^k=\sum_{n\geq 0}B(n,k)x^n$. Then the composita of the generating function $A(x)=B(x)-b(0)$ is equal to
\begin{equation}
A^{\Delta}(n,k)=\sum_{j=1}^k {k \choose j}B(n,j)(-1)^{k-j}b(0)^{k-j}.
\end{equation}
\end{theorem}

\begin{proof} Raising the generating function $A(x)$ to the power of $k$ and applying the binomial theorem, we obtain 
$$
A(x)^k=\left[B(x)-b(0)\right]^k=\sum_{j=0}^k {k \choose j}B(x)^j(-1)^{k-j}b(0)^{k-j}.
$$
From 
$$
[B(x)]^k=\sum_{n\geq 0}B(n,k)x^n,
$$
and
$B(x)^0=1$, we have
$$
B(n,0)=\begin{cases}
1, & \text{if $n=0$};\\
0, & \text{if $n>0$}.
\end{cases}
$$

Since $A(x)=\sum_{n>0} a(n)x^n$, we get
$$
A^{\Delta}(n,k)=\sum_{j=1}^k {k \choose j}B(n,j)(-1)^{k-j}b(0)^{k-j}.
$$
\end{proof}

\begin{theorem}
\label{Theorem_ak}
 Suppose $F(x)=\sum_{n>0} f(n)x^n$ is the generating function, $F^{\Delta}(n,k)$ is the composita of this generating function, and $\alpha$ is constant. Then for the generating function $A(x)=\alpha F(x)$  the composita is equal to 
\begin{equation}
A^{\Delta}(n,k)=\alpha^k F^{\Delta}(n,k).
\end{equation}
\end{theorem}
\begin{proof}
Using (\ref{GenComp}), we get
$$
[A(x)]^k=[\alpha F(x)]^k=\alpha^k[F(x)]^k=
$$
$$
=\sum_{n\geq  k} \alpha^k F^{\Delta}(n,k)x^n=\sum_{n\geq  k} A^{\Delta}(n,k)x^n.
$$
Therefore,
$$
A^{\Delta}(n,k)=\alpha^k F^{\Delta}(n,k).
$$
\end{proof}

\begin{theorem} \label{Theorem_an} 
Suppose $F(x)=\sum_{n>0} f(n)x^n$ is the generating function, $F^{\Delta}(n,k)$ is the composita of this generating function, and $\alpha$ is constant. Then for the generating function $A(x)=F(\alpha x)$ the composita is equal to
\begin{equation}
A^{\Delta}(n,k)=\alpha^n F^{\Delta}(n,k).
\end{equation}
\end{theorem}
\begin{proof}
Using (\ref{GenComp}), we get
$$
[A(x)]^k=[F(\alpha x)]^k=\sum_{n\geq  k} F^{\Delta}(n,k)(\alpha x)^n=
$$
$$
=\sum_{n\geq  k} \alpha^n F^{\Delta}(n,k)x^n=\sum_{n\geq  k} A^{\Delta}(n,k)x^n.
$$
Therefore,
$$
A^{\Delta}(n,k)=\alpha^n F^{\Delta}(n,k).
$$
\end{proof}

\begin{theorem} \label{Theorem_mult}
Suppose we have the generating function $F(x)=\sum_{n>0} f(n)x^n$, the composita of this generating function $F^{\Delta}(n,k)$; the following generating functions $B(x)=\sum_{n\geq 0} b(n)x^n$  and $[B(x)]^k=\sum_{n\geq 0}B(n,k)x^n$. 
Then for the generating function $A(x)=F(x)B(x)$  the composita is equal to 
\begin{equation}
\label{CompMult}
A^{\Delta}(n,k)=\sum_{i=k}^{n} F^{\Delta}(i,k)B(n-i,k).
\end{equation}
\end{theorem}

\begin{proof}
 Since $a(0)=f(0)b(0)=0$, the function $A(x)$ has the composita $A^{\Delta}(n,k)$. 
 
Using (\ref{GenComp}), we get
$$
 [A(x)]^k=[F(x)]^k[B(x)]^k.
$$ 
Then, from the rule of product of generating functions, we have
$$
A^{\Delta}(n,k)=\sum_{i=k}^{n} F^{\Delta}(i,k)B(n-i,k).
$$
\end{proof}

\begin{corollary}
If for the generating function $B(x)$ we have $b(0)=0$, then the formula (\ref{CompMult}) takes the form
\begin{equation}
A^{\Delta}(n,k)=\sum_{i=k}^{n-k} F^{\Delta}(i,k)B^{\Delta}(n-i,k).
\end{equation}
\end{corollary}

\begin{theorem} \label{Theorem_sum} 
Suppose we have the generating functions $F(x)=\sum_{n>0} f(n)x^n$, $G(x)=\sum_{n>0} g(n)x^n$, and their compositae $F^{\Delta}(n,k)$, $G^{\Delta}(n,k)$ respectively.
Then for the generating function $A(x)=F(x)+G(x)$  the composita is equal to
\begin{equation}
A^{\Delta}(n,k)=F^{\Delta}(n,k)+\sum_{j=1}^{k-1}{k\choose j}\sum_{i=j}^{n-k+j}F^{\Delta}(i,j)G^{\Delta}(n-i,k-j)+G^{\Delta}(n,k).
\end{equation}
\end{theorem}

\begin{proof}
Using (\ref{GenComp}) and the binomial theorem, we get
$$
[A(x)]^k=\sum_{j=0}^k{{k\choose j} [F(x)]^j[G(x)]^{k-j}}.
$$
Note that 
$$
[F(x)]^j=\sum_{n\geq j}F^{\Delta}(n,j)x^n,
$$
and
$$
[G(x)]^{k-j}=\sum_{n\geq k-j}G^{\Delta}(n,k-j)x^n.
$$
Then, from $F(x)^0=1$, $G(x)^0=1$ and the rule of product of generating functions, we have
$$
A^{\Delta}(n,k)=F^{\Delta}(n,k)+\sum_{j=1}^{k-1}{k\choose j}\sum_{i=j}^{n-k+j}F^{\Delta}(i,j)G^{\Delta}(n-i,k-j)+G^{\Delta}(n,k).
$$
\end{proof}

\begin{remark}
For the case $k=0$, we have $F(x)^0=1$. It is mean that 
\begin{equation}
F^{\Delta}(n,0)=\begin{cases}
1, & \text{if $n=0$};\\
0, & \text{if $n>0$}.
\end{cases}
\end{equation}
\end{remark}

\vspace{1.5cc}
\begin{center}
{\bf 3. COMPOSITAE OF GENERATING FUNCTIONS}
\end{center}

In this section we consider several examples of computation of compositae.

For derivation of a composita of the generating function $F(x)=\sum_{n>0} f(n)x^n$, we have to find coefficients of the generating function $F(x)^k$. As an example, in Table \ref{Tabl1_compositae} we present compositae of  several known generating functions \cite{Stanley_v2,ConcreteMath, Comtet}.

\begin{table}[h]
\begin{center}
\setlength\arrayrulewidth{1pt}
\renewcommand{\arraystretch}{1,3}
%\begin{tabular*}{0.9\linewidth}{|c @{\extracolsep{\fill}}c|}
\begin{tabular}{|ccc|ccc|}
\hline
&\textbf{Generating function $F(x)$}& & & \textbf{Composita $F^{\Delta}(n,k)$}&\\
\hline
&$x^m$ &&& $\delta_{n,mk}$, $m>0$&\\ \hline
&$\frac{bx}{1-ax}$ &&&$\binom{n-1}{k-1}a^{n-k}b^k$&\\  \hline
&$xe^x$ &&&$\frac{k^{n-k}}{(n-k)!}$&\\  \hline
&$\ln(1+x)$ &&& $\frac{k!}{n!}s(n,k)$& \\ \hline
&$e^x-1$  &&& $\frac{k!}{n!}S(n,k)$& \\ \hline
\end{tabular}
\caption{Examples of generating functions and their compositae}
\label{Tabl1_compositae}
\end{center}
\end{table}

Here  $\delta_{n,k}$ is the Kronecker delta, $s(n,k)$ and $S(n,k)$ stand for the Stirling numbers of the first kind and of the second kind, respectively (see \cite{Comtet,ConcreteMath}).

The Stirling numbers of the first kind $s(n,k)$ count the number of permutations of $n$ elements with $k$ disjoint cycles. The Stirling numbers of the first kind are defined by the following generating function
$$
\psi_k(x)=\sum_{n\geq k}  s(n,k)  \frac{x^n}{n!}=\frac{1}{k!}\ln^k(1+x).
$$

The Stirling numbers of the second kind $S(n,k)$ count the number of ways to partition a set of $n$ elements into $k$ nonempty subsets. A general formula for the Stirling numbers of the second kind is given as follows
$$S(n,k)=\frac{1}{k!}\sum_{j=0}^k(-1)^{k-j}\binom{k}{j}j^n.$$

The Stirling numbers of the second kind are defined by the following generating function
$$
\Phi_k(x)=\sum_{n\geq k} S(n,k) \frac{x^n}{n!}=\frac{1}{k!}(e^x-1)^k.
$$

\textbf{Compositae of polynomials}

Let us obtain compositae for polynomials. 
First we obtain the composita of the generating function $F(x)=ax+bx^2$. 
Raising this generating function to the power of $k$ and applying the binomial theorem, we get

$$[F(x,a,b)]^k=x^k(a+bx)^k=x^k\sum_{m=0}^k \binom{k}{m} a^{k-m}b^mx^m. 
$$

Substituting $n$ for $m+k$, we get the following expression
$$[F(x,a,b)]^k=\sum_{n=k}^{2k} \binom{k}{n-k}a^{2k-n}b^{n-k}x^n=\sum_{n=k}^{2k}F^{\Delta}(n,k,a,b)x^n. 
$$
Therefore, the composita is 
\begin{equation}
\label{Gnkab}
F^{\Delta}(n,k,a,b)=\binom{k}{n-k}a^{2k-n}b^{n-k}. 
\end{equation}

Next
we obtain the composita of the generating function $F(x)=ax+bx^2+cx^3$.
For this purpose, we write the generating function as the sum of the functions $F_1(x)=ax$ and $F_2(x)=x(bx+cx^2)$.

The composita of the generating function $F_1(x)=ax$, according to Theorem \ref{Theorem_ak}, is equal to 
$$F_1^{\Delta}(n,k,a)=a^k\delta_{n,k}.$$

Using Theorem \ref{Theorem_shift} and the formula (\ref{Gnkab}), the composita of the generating function  $F_2(x)$ is equal to
$$
F_2^{\Delta}(n,k,b,c)=\binom{k}{n-2k}b^{3k-n}c^{n-2k}. 
$$

Using Theorem \ref{Theorem_sum}, we obtain
$$
F^{\Delta}(n,k,a,b,c)=\sum_{j=0}^k{k\choose j}\sum_{i=j}^{n-k+j}F_1^{\Delta}(i,j,b,c)\delta_{n-i,k-j}a^{k-j}.
$$
Since
$$
\delta_{n-i,k-j}=\begin{cases}
 1,      & \text{if } n-i=k-j; \\
   0,        & \text{otherwise},
\end{cases}
$$
the composita of $F(x)=ax+bx^2+cx^3$ is
\begin{equation}
F^{\Delta}(n,k,a,b,c)=\sum_{j=0}^k{k\choose j}{j \choose n-k-j}a^{k-j}b^{2j+k-n}c^{n-k-j}.
\end{equation}

With the above theorems (Section 3), we can obtain compositae for different polynomials. Some examples are presented in Table \ref{Tabl2_compositae}.

\begin{table}[h]
\begin{center}
\setlength\arrayrulewidth{1pt}
\renewcommand{\arraystretch}{1,3}
%\begin{tabular*}{0.9\linewidth}{|c @{\extracolsep{\fill}}c|}
\begin{tabular}{|c|c|}
\hline
\textbf{Generating function $F(x)$}&  \textbf{Composita $F^{\Delta}(n,k)$}\\
\hline
$ax+bx^2$ & ${k \choose n-k} a^{2k-n}b^{n-k}$\\ \hline
$ax+bx^2+cx^3$ &$\sum\limits_{j=0}^k{k\choose j}{j \choose n-k-j}a^{k-j}b^{2j+k-n}c^{n-k-j}$\\  \hline
$ax+cx^3$ &$\frac{1+(-1)^{n-k}}{2}{k \choose \frac{n-k}{2}}a^\frac{3k-n}{2}c^\frac{n-k}{2}$ \\  \hline
$ax+bx^2+dx^4$ & $\sum\limits_{j=\left \lfloor {{4\,k-n}\over{3}} \right \rfloor}^{k}{a^{4\,k-n-2j}\,b^{n-4k+3j}d^{k-j}{{j}\choose{n-4k+3j}}\,{{k
 }\choose{j}}}$ \\ \hline
\end{tabular}
\caption{Compositae of polynomials}
\label{Tabl2_compositae}
\end{center}
\end{table}

\textbf{Compositae of trigonometric functions}

For computation of compositae of trigonometric functions, we use the Euler identity $e^{ix}=\cos(x)+i\sin(x)$.

Let us obtain the composita of the generating function $F(x)=\sin(x)$. 

Using the expression
$$
\sin(x)=\frac{e^{ix}-e^{-ix}}{2i},
$$
we obtain $\sin(x)^k$
\begin{multline*}
\sin(x)^k=\frac{1}{2^ki^k}\sum_{m=0}^k {k \choose m} e^{imx}e^{-i(k-m)x}(-1)^{k-m}=\frac{1}{2^ki^k}\sum_{m=0}^k {k \choose m}\times \\ \times e^{i(2m-k)x}(-1)^{k-m}=\sum_{n\geq k}\frac{1}{2^k}i^{n-k}\sum_{m=0}^k {k \choose m} \frac{(2m-k)^n}{n!}(-1)^{k-m}x^n.
\end{multline*}
Then the composita is equal to
$$
\frac{1}{2^k}i^{n-k}\sum_{m=0}^k {k \choose m} \frac{(2m-k)^n}{n!}(-1)^{k-m}.
$$

Since $n-k$ is an even number and the function is symmetric with respect to $k$, we obtain the composita of the generating function $F(x)=\sin (x)$
\begin{equation}
F^{\Delta}(n,k)=\begin{cases}
\frac{1}{2^{k-1}n!}\sum_{m=0}^{\lfloor \frac{k}{2} \rfloor} {k \choose m} (2m-k)^n(-1)^{\frac{n+k}{2}-m}, & \text{if $n-k$ is even};\\
0, & \text{if $n-k$ is odd}.
\end{cases}
\end{equation}

With the above theorems (Section 3), we can obtain compositae for different trigonometric and hyperbolic functions. Some examples are presented in Table \ref{Tabl3_compositae}.

\begin{table}[h]
\begin{center}
\setlength\arrayrulewidth{1pt}
\renewcommand{\arraystretch}{1,3}
%\begin{tabular*}{0.9\linewidth}{|c @{\extracolsep{\fill}}c|}
\begin{tabular}{|c|c|}
\hline
\textbf{Generating function $F(x)$}&  \textbf{Composita $F^{\Delta}(n,k)$}\\
\hline
$\sin(x)$ & $\frac{1+(-1)^{n-k}}{2^{k}n!}\sum\limits_{m=0}^{\lfloor \frac{k}{2} \rfloor} {k \choose m} (2m-k)^n(-1)^{\frac{n+k}{2}-m}
$\\ \hline
$x\cos(x)$ &$\frac{1+(-1)^{n-k}}{2^{k}(n-k)!}(-1)^{\frac{n-k}{2}}\sum\limits_{j=0}^{ \lfloor\frac{k-1}{2} \rfloor} {k \choose j}(k-2j)^{n-k}$\\  \hline
$\tan(x)$ &$
\frac{1+(-1)^{n-k}}{n!}\sum\limits_{j=k}^n 2^{n-j-1}\left\{{n \atop j}\right\}j!(-1)^{\frac{n+k}{2}+j}{j-1 \choose k-1}
$ \\  \hline
$\arctan(x)$ & $
\frac{\left((-1)^{\frac{3n+k}{2}}+(-1)^\frac{n-k}{2}\right)k!}{2^{k+1}}\sum\limits_{j=k}^n \frac{2^j}{j!}{n-1 \choose j-1}s(j,k)
$ \\ \hline
$\sinh(x)$ & $\frac{1}{2^k}\sum\limits_{j=0}^k (-1)^j{k\choose j} \frac{(k-2j)^n}{n!}$ \\ \hline
$x\cosh(x)$ & $\frac{1}{2^k}\sum\limits_{j=0}^k {k\choose j} \frac{(k-2j)^{n-k}}{(n-k)!}$ \\ \hline
\end{tabular}
\caption{Compositae of trigonometric and hyperbolic functions}
\label{Tabl3_compositae}
\end{center}
\end{table}

\vspace{1.5cc}
\begin{center}
{\bf 4. COMPOSITION OF GENERATING FUNCTIONS AND ITS COMPOSITAE}
\end{center}

Let us consider the application of compositae for computation of compositions of ordinary generating functions. For this purpose, we prove the following theorem.
\begin{theorem}
\label{Theorem_composition}
Suppose we have the generating function $F(x)=\sum_{n>0} f(n)x^n$, the composita of this generating function $F^{\Delta}(n,k)$, and the generating function $R(x)=\sum_{n\geq 0} r(n)x^n$. 
Then for the composition of generating functions $ A(x)=R(F(x))$ the following condition holds
\begin{equation}
\label{composition}
a(n)=\begin{cases}
r(0), & \text{if $n=0$};\\
\sum_{k=1}^n F^{\Delta}(n,k)r(k), & \text{if $n>0$},
\end{cases}
\end{equation}
where $A(x)=\sum_{n\geq 0} a(n)x^n.$
\end{theorem}

\begin{proof}
For computation $A(x)=R(F(x))$ we can write   
$$
A(x)=R(F(x))=\sum_{k\geq 0} r(k)F(x)^k.
$$
Replacing $F(x)^k$ by $\sum_{n\geq k}F^{\Delta}(n,k)x^n$ and considering that $F(x)^0=1$, we get 
$$
\begin{array}{llllll}
A(x)=r(0)+\\
&+r(1)F(1,1)x&+r(1)F(2,1)x^2&+\ldots+&r(1)F(n,1)x^n&+\cdots\\
 &          & +r(2)F(2,2)x^2&+\ldots+&r(2)F(n,2)x^n&+\cdots\\
  &         & & \cdots & \\
   &        & &+        & r(n)F(n,n)x^n&+\cdots\\
    &       & &  & &+\cdots
\end{array}
$$
Summing the coefficients of equal powers of $x^n$, we obtain the desired formula
$$
a(0)=r(0),\qquad n=0;
$$
$$
a(n)=\sum_{k=1}^nF^{\Delta}(n,k)r(k), \qquad n>0.
$$
\end{proof}

Further, for the composition $A(x)=R(F(x))$ the condition $a(0)=r(0)$ is implied.

\begin{example}
Let us obtain an expression of coefficients of the generating function 
$$A(x)=\frac{1}{1-ax-bx^2-cx^3},$$
 where $a,b,c\neq 0$.
 
Represent $A(x)$ as a composition of generating functions $A(x)=R(F(x))$,
where $F(x)=ax+bx^2+cx^3$ and $R(x)=\frac{1}{1-x}$.

According to Table \ref{Tabl2_compositae}, the composita of $F(x)=ax+bx^2+cx^3$ is
$$
\sum\limits_{j=0}^k{k\choose j}{j \choose n-k-j}a^{k-j}b^{2j+k-n}c^{n-k-j}.
$$

Using Theorem \ref{Theorem_composition}, we obtain the expression  of coefficients of $A(x)$ 
$$
a(n)=\sum_{k=1}^n \sum\limits_{j=0}^k{k\choose j}{j \choose n-k-j}a^{k-j}b^{2j+k-n}c^{n-k-j}.
$$
\end{example}

\begin{example}
Let us consider the generating function $A(x)=e^{\sinh(x)}$. 

Using  the composita of $F(x)=\sinh(x)$ (see Table \ref{Tabl3_compositae}) and Theorem \ref{Theorem_composition}, we obtain the expression  of coefficients of $A(x)$ 
$$
a(n)=\sum\limits_{k=1}^n\frac{1}{2^k}\sum\limits_{j=0}^k (-1)^j{k\choose j} \frac{(k-2j)^n}{n!}\frac{1}{k!}.
$$
\end{example}

\begin{theorem}\label{CompozitProduct}  
Suppose we have the generating functions $F(x)=\sum_{n>0}f(n)x^n$, $G(x)=\sum_{n>0} g(n)x^n$, and their compositae $F^{\Delta}(n,k)$, $G^{\Delta}(n,k)$ respectively.
Then for the composition of generating functions $A(x)=G\left(F(x)\right)$  the composita is equal to 
\begin{equation}
\label{compCompositon}
A^{\Delta}(n,k)=\sum_{m=k}^n F^{\Delta}(n,m)G^{\Delta}(m,k).
\end{equation}
\end{theorem}

\begin{proof}
Using the formula (\ref{GenComp}), we have
$$
[A(x)]^k=[G(F(x)]^k=\sum_{n\geq  k} A^{\Delta}(n,k)x^n.
$$

The function of coefficients of the generating function $[G(x)]^k$ is the composita $G^{\Delta}(n,k)$
$$
[G(x)]^k=\sum_{n\geq  k} G^{\Delta}(n,k)x^n.
$$

Then, according to Theorem \ref{Theorem_composition}, we get
$$
[G(F(x)]^k=\sum_{n\geq k}\sum_{m=1}^n F^{\Delta}(n,m)G^{\Delta}(m,k).
$$

Since 
$$
G^{\Delta}(m,k)=0, \qquad \text{if $m<k$},
$$
we obtain the composita of the composition of generating functions $A(x)=G\left(F(x)\right)$

$$
A^{\Delta}(n,k)=\sum_{m=k}^n F^{\Delta}(n,m)G^{\Delta}(m,k).
$$
\end{proof}

%\vspace{1.5cc}
\begin{center}
{\bf 5. COMPOSITAE OF RECIPROCAL GENERATING FUNCTIONS}
\end{center}

First we consider the notion of \textit{reciprocal generating functions} \cite{Wilf}.
\begin{definition}
\label{Recipro1}
Reciprocal generating functions are functions that satisfy the condition
$$
H(x)B(x)=1. 
$$
\end{definition}

\begin{remark}
If we have the reciprocal generating functions $H(x)=\sum_{n\geq 0} h(n)x^n$ and $B(x)=\sum_{n\geq 0} B(n)x^n$ such that $H(x)B(x)=1$, then by the composita of the reciprocal generating function of $B(x)$ we mean the composita of $xH(x)=\frac{x}{B(x)}$.
\end{remark}

In the following theorem we give the formula of the composita of a reciprocal generating function.
\begin{theorem}
\label{recip_Theorem}
Suppose  $H(x)=\sum_{n\geq 0}h(n)x^n$ is a generating function, $B(x)=\sum_{n\geq 0}b(n)x^n$ is the reciprocal generating function of $H(x)$, and $B_x^{\Delta}(n,k)$ is the composita of $xB(x)$. Then the composita of the generating function $xH(x)$ is equal to

\begin{equation}
H_x^{\Delta}(n,k)=
\begin{cases}
\frac{1}{B_x^{\Delta}(1,1)^{k}}, & \text{if $n=k$};\\
\sum\limits_{m=1}^{n-k}\binom{m+k-1}{k-1}
 \sum\limits_{j=1}^{m}\frac{(-1)^{j}}{B_x^{\Delta}(1,1)^{j+k}}\binom{m}{j}B_x^{\Delta}(n-k+j,j), &
 \text{if $n>k$}.
\end{cases}
\label{Reciprocal}
\end{equation}
\end{theorem}

\begin{proof} By Definition \ref{Recipro1}, we get  
$$
xH(x)=\frac{x}{b(0)+B(x)-b(0)}.
$$

Raising this generating function to the power of $k$, we obtain

$$
[xH(x)]^k=\left[\frac{x}{b(0)+B(x)-b(0)}\right]^k=\left[\frac{1}{b(0)}\frac{x}{1+\frac{1}{b(0)}(B(x)-b_0)}\right]^k.
$$

Using Corollary \ref{coeffB(n,k)}, Theorem \ref{comp-b} and Theorem \ref{Theorem_ak}, we obtain the composita of $F(x)=\frac{1}{b_0}(B(x)-b_0)$
$$
F^{\Delta}(n,k)=\sum_{j=1}^{k}{b(0)^{-j}\,(-1)^{k
 -j}\,{{k}\choose{j}}\,B_x^{\Delta}(n+j,j)}.
$$

The expression of coefficients of the generating function $R(x)=\left[\frac{1}{b(0)}\frac{1}{1+x}\right]^k$ is equal to
$$
R(n,k)=\frac{1}{b(0)^k}{n+k-1 \choose k-1}(-1)^{n}.
$$

Then, according to Theorem \ref{Theorem_composition}, we get
$$
H(n,k)=\begin{cases}
\frac{1}{b(0)^k}, & \text{if $n=0$};\\
\sum_{m=1}^n{m+k-1 \choose k-1}\sum_{j=1}^{m}\frac{1}{b(0)^{k+j}}(-1)^{j}\,{{m}\choose{j}}\,B^{\Delta}(n+j,j), & \text{if $n>0$}.
\end{cases}
$$

Therefore, from Corollary \ref{compB(n,k)} and $b(0)=B_x^{\Delta}(1,1)$, we obtain the composita of the reciprocal generating function of $B(x)$
$$
H_x^{\Delta}(n,k)=
\begin{cases}
\frac{1}{B_x^{\Delta}(1,1)^{k}}, & \text{if $n=k$};\\
\sum\limits_{m=1}^{n-k}\binom{m+k-1}{k-1}
 \sum\limits_{j=1}^{m}\frac{(-1)^{j}}{B_x^{\Delta}(1,1)^{j+k}}\binom{m}{j}B_x^{\Delta}(n-k+j,j), &
 \text{if $n>k$}.
\end{cases}
$$
\end{proof}

For applications of Theorem \ref{recip_Theorem} we give some examples.
\begin{example} Let us find a composita of the generating function $F(x)=x^2\csc(x)$. For this purpose, we write
$$
F(x)=x^2\csc(x)=\frac{x}{\frac{\sin(x)}{x}},
$$
or
$$
\frac{F(x)}{x}\frac{\sin(x)}{x}=1.
$$

According to Table \ref{Tabl3_compositae}, the composita of $\sin(x)$ is
$$
\frac{1+(-1)^{n-k}}{2^{k}n!}\sum_{m=0}^{\lfloor \frac{k}{2} \rfloor} {k \choose m} (2m-k)^n(-1)^{\frac{n+k}{2}-m}.
$$

Then, using Theorem \ref{recip_Theorem}, we obtain the composita of $F(x)$
$$
F^{\Delta}(n,k)=\begin{cases}
1, & \text{if $n=k$};\\
\sum\limits_{m=1}^{n-k}\binom{m+k-1}{k-1}
 \sum\limits_{j=1}^{m}\binom{m}{j}
 \frac{1+(-1)^{n-k}}{2^{j}(n-k+j)!}\sum\limits_{i=0}^{\lfloor \frac{j}{2} \rfloor} {j \choose i} (2i-j)^{n-k+j}(-1)^{\frac{n-k}{2}-i}, & \text{if $n>k$}.
\end{cases}
$$
\end{example}

\begin{example}
Let us find a composita of the generating function $F(x)=x\,H(x)$, where $H(x)=\sum_{n\geq 0}h(n)x^n$ is the generating function for Bernoulli numbers
$$
H(x)=\frac{x}{e^x-1}.
$$

According to Table \ref{Tabl1_compositae}, the composita of the generating function $e^x-1$ is equal to
$$
\frac{k!}{n!}S(n,k).
$$

The generating function $H(x)$ is the reciprocal generating function of $\frac{e^x-1}{x}$.
Then, using Theorem \ref{recip_Theorem}, we obtain the composita of the generating function $F(x)=\frac{x^2}{e^x-1}$
$$
F^{\Delta}(n,k)=\sum\limits_{m=0}^{n-k}\binom{m+k-1}{k-1}
 \sum\limits_{j=0}^{m}(-1)^{j}\binom{m}{j}\frac{j!}{(n-k+j)!}S(n-k+j,j).
$$
\end{example}

\vspace{1.5cc}
\begin{center}
{\bf 6. FUNCTIONAL EQUATION $B(x)=H(xB(x)^m)$}
\end{center}

First we consider a solution of the functional equation
\begin{equation}\label{Lagrange}
A(x)=xH(A(x)),
\end{equation}
where $A(x)$ and $H(x)$ are generating functions such that $H(x)=\sum_{n\geq 0}h(n)x^n$ and $A(x)=\sum_{n> 0}a(n)x^n$.

In the following lemma we give the Lagrange inversion formula, which was proved by Stanley \cite{Stanley_v2}.

\begin{lemma}[The Lagrange inversion formula]
\label{Lagrange formula}
Suppose $H(x)=\sum_{n\geq 0}h(n)x^n$ with $h(0)\neq 0$, and let $A(x)$ be defined by
\begin{equation}
\label{funcEq}
A(x)=xH(A(x)).
\end{equation}
Then
\begin{equation}
n[x^n]A(x)^k=k[x^{n-k}]H(x)^n,
\end{equation}
where $[x^n]A(x)^k$ is the coefficient of $x^n$ in  $A(x)^k$ and  $[x^{n-k}]H(x)^n$ is the coefficient of $x^{n-k}$ in  $H(x)^n$.
\end{lemma}

By using the above Lemma~\ref{Lagrange formula}, we now give the following theorem.

\begin{theorem}
\label{thm1} 
Suppose  $H(x)=\sum_{n\geq 0}h(n)x^n$ is a generating function, where $h(0)\neq 0$, $H_x^{\Delta}(n,k)$ is the composita of the generating function $xH(x)$, and $A(x)=\sum_{n>0}a(n)x^n$ is the generating function, which is obtained from the functional equation $A(x)=xH(A(x))$. Then the following condition holds true
\begin{equation}\label{LagrangeCompozita}
A^{\Delta}(n,k)=\frac{k}{n}H_x^{\Delta}(2n-k,n).
\end{equation}
\end{theorem}
\begin{proof} According to Lemma~\ref{Lagrange formula}, for the solution of the functional equation $A(x)=xH(A(x))$, we can write
$$
n[x^n]A(x)^k=k[x^{n-k}]H(x)^n.
$$

In the left-hand side, there is the composita of the generating function $A(x)$ multiplied by $n$:
$$
n[x^n]A(x)^k=n\,A^{\Delta}(n,k).
$$

We know that
$$
\left( xH(x)\right)^k=\sum_{n\geq k}H_x^{\Delta}(n,k)x^n.
$$
Then
$$
\left( H(x)\right)^k=\sum_{n\geq k}H_x^{\Delta}(n,k)x^{n-k}.
$$
If we replace $n-k$ by $m$, we obtain the following expression
$$
\left( H(x)\right)^k=\sum_{m\geq 0}H_x^{\Delta}(m+k,k)x^{m}.
$$
Substituting $n$ for $k$ and $n-k$ for $m$, we get
$$
[x^{n-k}]H(x)^n=H_x^{\Delta}(2n-k,n).
$$
Therefore, we get
$$
A^{\Delta}(n,k)=\frac{k}{n}H_x^{\Delta}(2n-k,n).
$$
\end{proof}

According to the above Theorem, for solutions of the functional equation $A(x)=xH(A(x))$, we can use the following expression
$$
[A(x)]^k=\sum_{n\geq k}A^{\Delta}(n,k)x^n=\sum_{n\geq k} \frac{k}{n}H_x^{\Delta}(2n-k,n)x^n,
$$
where $H_x^{\Delta}(n,k)$ is the composita of the generating function $xH(x)$.
Therefore,
\begin{equation}
A(x)=\sum_{n\geq 1} \frac{1}{n}H_x^{\Delta}(2n-1,n)x^n.
\label{solvEq1}
\end{equation}

Since the composita is a unique of the given generating function, the formula (\ref{LagrangeCompozita}) provides a solution of the inverse equation $A(x)=xH(A(x))$, when $A(x)$ is known and $H(x)$ is unknown. Hence,
$$
H_x^{\Delta}(n,k)=\frac{k}{2k-n}A^{\Delta}(k,2k-n).
$$
It should be noted that for $n=k$,
$$
H_x^{\Delta}(n,n)=A^{\Delta}(n,n).
$$

Next we give some examples of functional equations.

\begin{example} 
Let us find the generating function $A(x)=\sum_{n>0}a(n)x^n$, which is defined by the functional equation
$$
A(x)=x+xA(x)+xA(x)^2+2xA(x)^3. 
$$
The generating function $xH(x)$ has the form
$$
xH(x)=x+x^2+x^3+2x^4.
$$
The composita of $xH(x)$  is
$$
H_x^{\Delta}(n,k)=\sum\limits_{j=0}^{k}{{{k}\choose{j}}\,\sum\limits_{i=j}^{n-k+j}{2^{n-3(k-j)-i}\,{{j
 }\choose{i-j}}{{k-j}\choose{n-3(k-j)-i}}}}.
$$
According to the formula (\ref{solvEq1}), the coefficients of $A(x)$ is
$$a(n)=\frac{1}{n}G^{\Delta}(2n-1,n).$$

Therefore, we get 
$$
a(n)=\frac{1}{n}{{\sum\limits_{j=0}^{n}{{{n}\choose{j}}\,\sum\limits_{i=j}^{n+j-1}{{{j}\choose{i-
 j}}\,2^{-n+3\,j-i-1}\,{{n-j}\choose{-n+3\,j-i-1}}}}}}.
$$
\end{example}

\begin{example} 
Let us find coefficients of the generating function $B(x)=\sum_{n\geq 0}a(n)x^n$, which is defined by the functional equation (see A064641 \cite{oeis})
$$
B(x)=\frac{1-x-\sqrt{1-6x-3x^2}}{ 2x(1+x)}
$$

Next we introduce the following generating function $A(x)=xB(x)$.
Considering the functional equation, we can take notice that
$$A(x) =x\frac{1+A(x) +A(x) ^2}{1-A(x)}.$$

Then we get the following functional equation
$$
A(x) = xH(A(x)),
$$
where $H(x) = \frac{1+x+x^2}{1-x}$.

Now we obtain the composita of $xH(x)$.
According to Table \ref{Tabl2_compositae}, the composita of $F(x)= x+x^2+x^3$ is
$$
F^{\Delta}(n,k)=\sum_{j=0}^k {j \choose n-k-j}{k \choose j}. 
$$

The expression of coefficients of the generating function $[R(x)]^k=\left(\frac{1}{1-x}\right)^k$ is equal to
$${n+k-1 \choose k-1}.$$

Then, using Theorem \ref{Theorem_mult}, we obtain the composita of $xH(x)$
$$
H_x^{\Delta}(n,k)=\sum_{i=0}^{n-k}{k+i-1 \choose k-1}\sum_{j=0}^{k}{j \choose n-k-j-i}{k \choose j}.
$$

Hence, using Theorem \ref{thm1}, we get the composita of $A(x)$
$$
A^{\Delta}(n,k)=\frac{k}{n}H_x^{\Delta}(2n-k,n)=\frac{k}{n}\sum_{i=0}^{n-k}{n+i-1 \choose n-1}\sum_{j=0}^{n}{j \choose n-k-j-i}{n \choose j}.
$$

According to Theorem \ref{coeffB(n,k)}, the coefficients of $B(x)$ is
$$b(n)=A^{\Delta}(n+1,1).$$

Therefore, we obtain
$$
b(n)=\frac{1}{n+1}\sum_{i=0}^{n}{n+i \choose n}\sum_{j=0}^{n+1}{j \choose n-j-i}{n+1 \choose j}.
$$
\end{example}

Next we generalize the case $A(x)=xH(A(x))$. 

Replacing $A(x)$ by $xB(x)$ in the functional equation (\ref{funcEq}), we get
\begin{equation}
B(x)=H(xB(x)).
\label{B(x)eq}
\end{equation}

Let us introduce the following definitions.
\begin{definition} The left composita of the generating function $B(x)$ in the functional equation (\ref{B(x)eq}) is the composita
$$
H_x^{\Delta}(n,k)=\frac{k}{2k-n}B_x^{\Delta}(k,2k-n),
$$
where $H_x^{\Delta}(n,k)$ is the composita of the generating function $xB(x)$.
\end{definition}

\begin{definition} The right composita of the generating function $H(x)$ in the functional equation (\ref{B(x)eq}) is the composita
$$
B_x^{\Delta}(n,k)=\frac{k}{n}H_x^{\Delta}(2n-k,n),
$$
where $H_x^{\Delta}(n,k)$ is the composita of the generating function $xH(x)$.
\end{definition}
There exists the left composita for every left composita and there exists the right composita for every right composita.

The formula (\ref{LagrangeCompozita}) can be generalized for the case generating function is the solution of a certain functional equation. 
Let us prove the following theorem.

\begin{theorem}
\label{Generalization}
Suppose  $H(x)=\sum_{n\geq 0}h(n)x^n$ and $B(x)=\sum_{n\geq 0}b(n)x^n$ are generating functions such that  
$B(x)=H(xB(x)^m)$, where $m\in \mathbb{N}$; $H_x^{\Delta}(n,k)$ and $B_x^{\Delta}(n,k)$ are the compositae of the generating functions $xH(x)$ and $xB(x)$, respectively. Then
\begin{equation}
B_x^{\Delta}(n,k)=\frac{k}{i_{m-1}}H_x^{\Delta}(i_m,i_{m-1}),
\end{equation}
where $i_m=(m+1)n-mk.
$
\end{theorem} 

\begin{proof}
For $m=0$, we have
$$
B(x)=H(xB(x)^0)=G(x), \qquad i_{m-1}=k, i_m=n.
$$
Then we obtain the identity
$$B_x^{\Delta}(n,k)=\frac{k}{k}H_x^{\Delta}(n,k).$$

For $m=1$, we have
$$
B_1(x)=G(xB_1(x)), \qquad i_{m-1}=n, i_m=2n-k.
$$
Then we obtain
$$B_{x, 1}^{\Delta}(n,k)=\frac{k}{n}H_x^{\Delta}(2n-k,n)$$
that satisfy Theorem \ref{LagrangeCompozita}.

By induction, we put that for $m$ the solution of the equation
\begin{equation}\label{SolverInitT1}
B_m(x)=H(xB_m(x)^m)
\end{equation}
is
$$
B_{x, m}^{\Delta}(n,k)=\frac{k}{i_{m-1}}H_x^{\Delta}(i_m,i_{m-1}).
$$

Then we find the solution for $m+1$
$$
B_{m+1}(x)=H(xB_{m+1}(x)^{m+1}).
$$
For this purpose, we consider the following functional equation
$$
B_{m+1}(x)=B_m(xB_{m+1}(x)).
$$

Instead of $B_m(x)$ we substitute the right hand-side of (\ref{SolverInitT1})
$$
B_{m+1}(x)=H(xB_{m+1}(x)[B_m(xB_{m+1}(x)]^m);
$$
from whence it follows that
$$
B_{m+1}(x)=H(xB_{m+1}(x)^{m+1}).
$$

We note that $B_{x, m+1}^{\Delta}(n,k)$ is the right composita of $B_m(x)$
$$
B_{x, m+1}^{\Delta}(n,k)=\frac{k}{n}B_{x, m}^{\Delta}(2n-k,n),
$$
where $B_{x, m}^{\Delta}(n,k)$ is the composita of the generating function $xB_m(x)$.

Then
$$
B_{x, m+1}^{\Delta}(n,k)=\frac{k}{(m+1)n-mk}H_x^{\Delta}((m+2)n-(m+1)k, (m+1)n-mk).
$$

Therefore, for  the functional equation
$$
B_{m+1}(x)=B_m(xB_{m+1}(x))
$$
we obtain the required condition
$$
B_{x, m+1}^{\Delta}(n,k)=\frac{k}{i_{m}}H_x^{\Delta}(i_{m+1},i_{m}),
$$
where $i_m=(m+1)n-mk$.
\end{proof}

In Table \ref{tablepascal} we present a sequence of functional equations for the generating function $H(x)=1+x$.

\begin{table}
\begin{center}\begin{tabular}{|c|c|c|c|} \hline
Equation & Function $B(x)$ & Composita $B_{x}^{\Delta}(n,k)$ & OEIS\\  \hline
$B(x)=1+xB^0(x)$ & $1+x$ & ${k \choose n-k}$& \\ \hline
$B(x)=1+xB^1(x)$ & $\frac{1}{1-x}$ & ${n-1 \choose k-1}$ & A000012\\ \hline
$B(x)=1+xB^2(x)$ & $\frac{1-\sqrt{1-4x}}{2x} $ &$\frac{k}{n}{2n-k-1 \choose n-1}$ & A000108 \\ \hline
$B(x)=1+xB^3(x)$ & & $\frac{k}{n}{3n-2k \choose n-k}$& A001764\\ \hline
\end{tabular}
\caption{Table of functional equations}\label{tablepascal}
\end{center}
\end{table}

\vspace{1.5cc}
\begin{center}
{\bf 7. CONCLUSION}
\end{center}
In this paper we introduce the concept of composita for ordinary generating functions and provide a number of applications. The proposed apparatus of compositae is applicable to solve the following problems: calculation of the composition of ordinary generating functions; finding expressions of reciprocal generating functions; finding expressions of inverse generating functions; finding solutions of functional equations; obtaining of expressions polynomials and etc.

\vspace{1cc}

Tomsk State University of Control Systems and Radioelectronics(TUSUR),\\
Tomsk, Russia 634050\\
email: KruchininDm@gmail.com (Dmitry Kruchinin);
kru@ie.tusur.ru (Vladimir Kruchinin).

{\small
\noindent

}
\vspace{2cc}


\begin{thebibliography}{1}

\bibitem{Stanley_v2} R.~P. Stanley.
\newblock {\em Enumerative combinatorics $2$}, volume~62 of { Cambridge
  Studies in Advanced Mathematics}.
\newblock Cambridge University Press, Cambridge, 1999.

\bibitem{ConcreteMath} R.~L. Graham, D.~E. Knuth, and O.~Patashnik, \emph{Concrete mathematics}, \newblock Addison-Wesley, Reading, MA, 1989.

\bibitem{Egor} G.~P. Egorichev, {\em Integral representation and the computation of combinatorial sums} , \newblock Amer. Math. Soc.  (1984)  

\bibitem{WJohnson}
W.~P. Johnson, {\em The Curious History of Faa du Bruno's Formula} //
\newblock The American Mathematical Monthly, vol. 109, 2002, pp. 217-234


\bibitem{Kru2010} V.~V. Kruchinin,
\newblock {\em Combinatorics of Compositions and its Applications}, \newblock V-Spektr, Tomsk, 2010. (in rus)

\bibitem{Shapiro} L. W. Shapiro, S. Getu, W.-J. Woan, and L. Woodson, {\em The Riordan group}, \newblock Discrete
Applied Math. 34 (1991), 229{339}.

\bibitem{Sprugnoli} R. Sprugnoli. {\em Riordan arrays and combinatorial sums}. Discrete Mathematics,
132:267{290, 1994}.

\bibitem{Comtet} L. Comtet, 
\newblock \textit{Advanced Combinatorics}, D. Reidel Publishing Company, 1974.

\bibitem{Flajolet} P. Flajolet and R. Sedgewick, \emph{Analytic Combinatorics}, Cambridge University Press, 2008.


\bibitem{Wilf} H.~S. Wilf {\em Generatingfunctionology} \newblock Academic Press, 1994.

\bibitem{oeis}
N. J. A. Sloane, \emph{The On-Line Encyclopedia of Integer Sequences},
http://oeis.org.
\end{thebibliography}
\end{document}